\documentclass{amsart}
\usepackage{amssymb,  amscd, times}

\def\CC{{\mathbb C}}

\def\PP{{\mathbb P}}
\def\QQ{{\mathbb Q}} 
\def\RR{{\mathbb R}}

\def\ZZ{{\mathbb Z}}

\def\G{\Gamma}

\def\bs{\backslash}
\newcommand{\eps}{\varepsilon}

\def\Ccal{{\mathcal C}}

\def\Hcal{{\mathcal H}} 
\def\Ical{{\mathcal I}} 
 
\def\Kcal{{\mathcal K}}
\def\Lcal{{\mathcal L}}

\def\Ocal{{\mathcal O}}

\def\Scal{{\mathcal S}}

\def\Qcal{{\mathcal Q}}

\def\Zcal{{\mathcal Z}}

\def\s{{\rm s}}
\def\ss{{\rm ss}}

\def\pt{{\scriptscriptstyle\bullet}}

\newcommand\aut{\operatorname{Aut}}

\newcommand\iso{\operatorname{Iso}}

\newcommand\pic{\operatorname{Pic}}
\newcommand\proj{\operatorname{Proj}}

\newcommand\res{\operatorname{Res}}
\newcommand\spec{\operatorname{Spec}}

\newcommand\sym{\operatorname{Sym}}

\newcommand\GL{\operatorname{GL}}
\newcommand\Or{\operatorname{O}}
\newcommand\SL{\operatorname{SL}}
\newcommand\PL{\operatorname{PGL}}

\newcommand\U{\operatorname{U}}

\newcommand\sign{\operatorname{sign}}

\newtheorem{theorem}{Theorem}[section]
\newtheorem{lemma}[theorem]{Lemma}
\newtheorem{proposition}[theorem]{Proposition}

\newtheorem{corollary}[theorem]{Corollary}

\theoremstyle{definition}
\newtheorem{definition}[theorem]{Definition}

\theoremstyle{remark} 
\newtheorem{remark}[theorem]{Remark}

\newtheorem{question}[theorem]{Question}

\title{Invariants of quartic plane curves as automorphic forms}
\author{Eduard Looijenga}
\address{Mathematisch Instituut\\
Universiteit Utrecht\\
P.O.~Box 80.010, NL-3508 TA Utrecht\\
Nederland}
\email{looijeng@math.uu.nl}
\begin{document}

\begin{abstract}
We identify the algebra of  regular functions  on the space of quartic polynomials in three complex variables invariant under $\SL (3,\CC)$ with an algebra of meromorphic automorphic forms
on the complex $6$-ball. We also discuss the underlying geometry.
\end{abstract}
\maketitle
\hfill{\emph{To Igor Dolgachev, for his 60th birthday}}
\section{Motivation and goal} 
One of the most beautiful mathematical gems of the 19th century is,
at least to my taste (but I expect Igor to agree), the theorem 
that says that the algebra of invariants of plane cubics  is the algebra of modular
forms. Its precise statement is as follows. Consider the  vector space of cubic homogeneous forms 
in three complex variables, which as a $\SL (3,\CC)$-representation is the third symmetric power of the dual of the 
defining representation $\CC^3$. Since $\SL (3,\CC)$ acts through 
its simple quotient $\PL (3,\CC)$, we prefer to regard this a representation of the latter.
Then the theorem I am refering to says that the algebra of  $\PL (3,\CC)$-invariant polynomials on that space is as a graded $\CC$-algebra isomorphic to the algebra of $\SL(2,\ZZ)$-modular forms  (recall that this
is the polynomial algebra generated by the Eisenstein series 
$E_4$ and $E_6$). This theorem is in a sense the optimal algebraic form of the
geometric property which says that a genus one curve is completely determined by its periods. Yet there is subtlety here which 
deserves to be explicated. Denote the vector space of 
homogeneous cubic forms on $\CC^3$ by $\Kcal$
and let $\Kcal^\circ\subset \Kcal$ be the open subset of those $F\in \Kcal$ 
whose zero set $C(F)$ in $\PP^2$ 
defines a smooth cubic curve (hence of genus one). Then any $F\in \Kcal^\circ$ 
also  determines a distinguished holomorphic  differential $\omega (F)$ on 
$C(F)$, which perhaps is best  described as an iterated residue:
\[
\omega(F)=\res_{C(F)}\res_{Z(F)} \frac{dZ_0\wedge dZ_1\wedge dZ_2}{F}.
\]
Here $Z(F)$ denotes the zero set of $F$ in $\CC^3$ (the affine cone over
$C(F)$) and $C(F)$ is regarded as
the locus where $Z(F)$ meets the hyperplane at infinity $\PP^2$. It is easy
to see  that this double residue is indeed a nonzero holomorphic differential on 
$C(F)$. The lattice of periods of $\omega (F)$ (a lattice in $\CC$)
only depends on the $\PL (3,\CC)$-orbit of $F$.
Now recall that the space of  lattices in $\CC$ is naturally the $\SL(2,\ZZ)$-orbit space of the space $\iso^+(\RR^2,\CC)$ of oriented $\RR$-isomorphisms $\RR^2\to\CC$  (where $\sigma\in \SL(2,\ZZ)$ takes $\zeta\in \iso^+(\RR^2,\CC)$ to $\zeta\sigma^{-1}$). 
So we have defined a map
\[
\PL(3,\CC)\bs \Kcal^\circ\to \SL(2,\ZZ)\bs \iso^+(\RR^2,\CC).
\]
It is easy to see that this map is in fact an isomorphism of (affine) algebraic varieties
(as for its injectivity, observe that the $\PL (3,\CC)$-stabilizer of 
a nonsingular cubic plane curve  is the group of automorphisms of that curve which preserve its degree three polarization). Both sides come with a natural $\CC^\times$-action: $\lambda\in\CC^\times$ acts on
$\Kcal^\circ$ by multiplication with $\lambda^{-1}$ and acts on 
$\iso^+(\RR^2,\CC)$ by composing with multiplication by $\lambda$ on $\CC$.
These actions descend to the orbit spaces and make the isomorphism above 
$\CC^\times$-equivariant. (These actions are no longer effective, for
$-1\in\CC^\times$ acts as the identity.)
The graded algebra $\CC [\Kcal]^{\PL(3,\CC)}$ can be understood as 
an algebra of regular functions on the left hand side and the graded algebra
of $\SL(2,\ZZ)$-modular forms $\CC[E_4,E_6]$ can be understood as 
an algebra of regular functions on right hand side. In either case the 
grading comes from the $\CC^\times$-action we just described. 
The cited theorem says that the displayed map identifies these  
graded $\CC$-algebras.
This has a geometrical consequence which goes somewhat 
beyond the observed  isomorphism: the proj construction on the left is 
interpreted by Geometric Invariant Theory: 
$\proj (\CC [\Kcal]^{\PL (3,\CC)})$  adds to $\PL(3,\CC)\bs \PP(\Kcal^\circ)$ 
a singleton which is represented by the unique closed strictly semistable 
$\PL(3,\CC)$-orbit in $\PP(\Kcal)$, namely, the orbit of   three nonconcurrent lines in $\PP^2$. 
The proj construction on the 
right compactifies the modular curve in question, that is, the $j$-line  
$\SL(2,\ZZ)\bs \PP\iso^+(\RR^2,\CC)$, in the standard manner: it is the
simplest instance of a Baily-Borel compactification.\\ 

To me this theorem serves as a model for the theory of period 
maps. It seems to tell  us that whenever we know such a map to be an 
open embedding, then we should try to express that fact on a (richer) 
algebraic level, amounting to the identification of two graded algebras of 
invariants, one with respect to a reductive  algebraic group, the other
with respect to a discrete group.
The ensueing identification of their 
proj's should give us then an additional piece of geometric information,
which in the end is a sophisticated way to understand 
the period map's boundary behavior. 
This is most likely to 
happen when that map takes values in  a ball quotient or a locally 
symmetric variety of type IV. If the period map is an isomorphism, then one
hopes for an exact analogue of the theorem above. For example, 
Allcock-Carlsen-Toledo \cite{act} have essentially established this 
for the case of cubic surfaces, where $\Kcal$ is replaced by the space of 
homogenenous cubic forms in 4 variables and the algebra of modular forms 
by the algebra of automorphic forms on a $4$-ball with respect to an arithmetic group. 

But if the period map is not
surjective, then some modifications  on the automorphic side are in order 
and it is precisely for this purpose that I developed the geometric theory 
of meromorphic automorphic forms in \cite{looij1} and \cite{looij2}. What I 
want to do here is to illustrate that theory in the (other) case that logically 
comes after our guiding example, namely that of quartic plane curves. 
To be precise, let $\Qcal$ stands for the vector space of \emph{quartic} 
(rather than cubic) homogeneous forms in three complex variables and regard this
space as a representation of $\SL (3,\CC)$, then
we shall interpret the algebra of invariants $\CC [\Qcal]^{\SL(3,\CC)}$ 
(or rather $\CC [\Qcal]^{\mu_4\times\SL(3,\CC)}$ with $\mu_4$ acting
on $\Qcal$ by scalar multiplication) as an 
algebra of (meromorphic) automorphic forms on the complex $6$-ball. 
This will be an algebra isomorphism which rescales the degrees by a factor three. 
We shall, of course, also interpret the proj construction on either side.

This example is not an isolated one. For instance, Allcock \cite{allcock} has 
determined the  semistable cubic threefolds and his results suggest that the 
situation is very much like the case of quartic curves. \\

It is pleasure to  dedicate this paper to my longtime friend Igor Dolgachev.
Igor and I share a passion for our field (which is not confined to algebraic geometry) and we have a similar mathematical taste, but  I only wish I had his extensive knowledge of the classical literature of our subject. I learned a lot from him. 

\section{The Baily-Borel compactification of a ball quotient}

Let $V$ be a complex vector space of finite dimension, endowed with a nondegenerate Hermitian form $h:V\times V\to \CC$ of signature $(1,n)$. So we can find a coordinate system $(z_0,\dots ,z_n)$  for $V$ on which $h$ takes the standard form:
$h(z,z)=|z_0|^2-|z_1|^2-\cdots -|z_n|^2$. This shows that the 
unitary group $\U(V)$ of $(V,h)$ is isomorphic to $\U (1,n)$.
Let us denote by $V_+$ the open subset of $z\in V$ with $h(z,z)>0$. 
This set is $\CC^\times$-invariant and hence defines an open subset
$\PP (V_+$) of $\PP (V)$. In terms of the above coordinates, $\PP (V_+)$ is
defined by $\sum_{i=1}^n |z_i/z_0|^2<1$, which shows that $\PP (V_+)$ is 
biholomorphic to the complex $n$-ball.  The group $\U(V)$ acts properly
and transitively on it; in fact, the stabilizer of a point is a maximal 
compact subgroup $\U(V)$, so that $\PP( V_+)$ can be understood as the
symmetric space of $\U(V)$. 

Suppose we are also given a discrete subgroup $\Gamma\subset \U(V)$ of
finite \emph{covolume} (which means that $\U(V)/\Gamma$ has finite 
$\U(V)$-invariant volume). Since  $\G$ is discrete, it acts properly 
discontinously on $\PP (V_+)$ (and hence also on $V_+$) so that
the formation of the orbit spaces $\G\bs \PP (V_+)$ and $\G\bs V_+$ 
takes place in the complex-analytic category (the orbit spaces
will be normal). However, the Baily-Borel theory tells us that
these spaces have the richer structure of quasi-projective variety.
In order to state a precise result, we recall that a \emph{$\G$-automorphic form
of degree $d\in\ZZ$} is in the present setting a $\G$-invariant holomorphic function
$f:V_+\to \CC$ which is homogeneous of degree $-d$: $f(\lambda z)=
\lambda^{-d}f(z)$ and which in case $n=1$ also obeys a growth condition
which we do not bother to specify. If we denote the space of such forms by $A_d^\G$, 
then it is clear that their direct sum $A_\pt^\G$ is a graded $\CC$-algebra.
The Baily-Borel theory has the following to say about it:

\begin{theorem}\label{thm:BB}
The graded $\CC$-algebra $A_\pt^\G$ has finitely many generators of positive
degree. This algebra separates the $\G$-orbits in $V_+$ so that we have
injective complex-analytic maps
\[
\G\bs V_+\to \spec (A_\pt^\G),\quad \G\bs \PP( V_+)\to \proj (A_\pt^\G).
\]
The images of these morphisms are Zariski open-dense so that
$\G\bs V_+$ resp.\ $\G\bs \PP (V_+)$ acquires the structure of a
quasi-affine resp.\ quasi-projective complex variety. Moreover,
we have natural (so-called \emph{Baily-Borel}) extensions $V_+\subset V_+^*$ and 
$\PP( V_+)\subset 
\PP (V_+^*)$ as topological spaces with $\G$-action (which we 
describe below) such that these morphisms extend to homeomorphisms 
\[
\G\bs V^*_+\cong \spec (A_\pt^\G),
\quad \G\bs \PP (V^*_+)\cong \proj (A_\pt^\G),
\]
if we endow the targets with their Hausdorff topology.
\end{theorem} 
 
In order to describe the extensions mentioned in this theorem, 
we introduce the following notion. Let $H\subset V$ be a linear hyperplane.
If $H$ is nondegenerate, then the $\U(V)$-stabilizer $\U(V)_H$ is simply
$U(H)\times U(H^\perp)$. Suppose now that $H$ is degenerate. Then 
$H^\perp\subset H$,  $H/H^\perp$ is negative definite and the evident 
homomorphism $\U(V)_H\to\U (H/H^\perp)\times \GL (H^\perp)$ is surjective with kernel a 
unipotent group (the unipotent radical
of $\U(V)_H$).  This unipotent group is in fact a Heisenberg group: 
an extension of the complex vector group 
$(H/H^\perp)\otimes \overline{H^\perp}$ 
by a real vector group of dimension one.

\begin{definition}
We call a linear hyperplane $H\subset V$ that is not negative 
definite \emph{$\G$-rational} if its $\G$-stabilizer $\G_H$ maps to
a subgroup of finite covolume in $\U( H)$.
\end{definition}

(Since the  $\U(V)$-stabilizer of a negative definite hyperplane 
is  compact,  this notion would be for such hyperplanes
without much interest.) Let $\Ical$ denote the collection of isotropic 
lines $I\subset V$ for which $I^\perp$ is 
$\G$-rational. The $\G$-rationality amounts here to requiring  that 
$\G_I$ meets the unipotent radical of  $\U(V)_I$ in a cocompact subgroup. 
According to the reduction theory for such forms, $\G$ has only finitely 
many orbits in the set $\Ical$. Now we can
describe the extensions mentioned in Theorem \ref{thm:BB} as a set:
\[
V^*_+=V_+\sqcup \left(\coprod_{I\in\Ical} V_+/I^\perp\right)\sqcup (V_+/V),\quad 
\PP (V^*_+)=\PP (V_+)\sqcup \coprod_{I\in\Ical}\PP ( V_+/I^\perp),
\]
where $V_+/I^\perp$  simply denotes the image of $V_+$ in $V/I^\perp$ and likewise 
in other cases. This notation may appear unnecessarily complicated, since  
$V_+/I^\perp= V/I^\perp-\{ 0\}$ is just a copy of  $\CC^\times$ and both $V_+/V$ 
and $\PP ( V_+/I^\perp)$ are even singletons. Yet it is convenient notation, not just 
because we need to be able to distinguish singletons by name, but also 
because it is typical for the general situation that the extension is 
obtained by adding quotients of the very space we are extending. Still we 
may observe that
$\PP (V^*_+)-\PP (V_+)$ can be identified with the discrete set $\Ical$ 
so that $\G\bs \PP (V^*_+) -\G\bs  \PP (V_+)$ is identified with 
the finite set $\G\bs \Ical$.

We will not define the topology of these extensions; suffices to say here 
that this topology induces the given one on
each stratum, and  that the elements of both $\G$ and $\CC^\times$ 
(acting by scalar multiplication) act as homeomorphisms. The space $V^*_+$ 
has the structure of a cone
with  the singleton $V_+/V$ as vertex and $\PP( V^*_+)$ as its base.

\section{$\G$-arrangements and associated compactifications}

We keep the situation of the previous section, but now assume that $n\ge 2$. 

\subsection*{$\G$-arrangements}
Suppose that $H\subset V$ is a 
$\G$-rational hyperplane  with $H^\perp$ a negative definite line. So $H$ has 
signature $(1,n-1)$, $H_+=H\cap V_+$ and $\G_H$ maps to a subgroup
of $U(H)$ of cofinite volume. That subgroup will be discrete and so the preceding applies: we  find a quasi-affine variety  $\G_H\bs H_+$ with its 
Baily-Borel extension $\G_H\bs H_+^*$ (a normal affine variety). 
The natural map $\G_H\bs H_+\to \G\bs V_+$ is evidently complex-analytic, 
but since $\G$-automorphic forms restrict to $\G_H$-automorphic forms
it is in fact an algebraic morphism. For the same reason,
this map has a Baily-Borel extension $\spec(A_\pt^{\G_H})\to \spec(A_\pt^\G)$ 
in the algebraic category. The underlying map $\G_H\bs H^*_+\to \G\bs V^*_+$
is defined in an evident manner. In particular, the preimage of the vertex is the vertex and so this Baily-Borel extension is finite. It is also  birational onto its image 
and hence this map is a normalization of that image. Since the preimage of $\G\bs V_+$ in  $\G_H\bs H^*_+$ is  
$\G_H\bs H_+$, it follows that $\G_H\bs H_+\to \G\bs V_+$ is a finite morphism
as well. The preimage in $V_+$ of the image of  
$\G_H\bs H_+\to \G\bs V_+$ is of course the union of the $\G$-translates of
$H_+$. It follows that these translates form a  collection of hyperplane
sections of $V_+$ that is locally finite. This remains true if instead
of a single $H$ we take a finite number of them.
This leads up to the following

\begin{definition}
A \emph{$\G$-arrangement} is a collection $\Hcal$ of $\G$-rational hyperplanes  
with negative definite orthogonal complement, which, when viewed
as a subset of the appropriate Grassmannian, is a union of finitely many $\G$-orbits.
\end{definition}

For the remainder of this section we fix a $\G$-arrangement $\Hcal$.

\subsection*{Meromorphic automorphic forms}
 The preceding discussion shows that
the collection $(H_+)_{H\in\Hcal}$ is locally finite on $V_+$ so that its
union $D_\Hcal:=\cup_{H\in\Hcal} H_+$ is closed in $V_+$. Moreover,
$D_\Hcal$ is the preimage of a hypersurface in $\G\bs V_+$. We shall denote
$V_+-D_\Hcal$ simply by $V_{\Hcal}$. Any open subset of $V$ of this form will be
refered to as a \emph{$\G$-arrangement complement}. It is clear that 
$\G\bs V_{\Hcal}$ is the complement of a hypersurface in $\G\bs V_+$. It is our goal 
to state a generalization of Theorem \ref{thm:BB} for $\G\bs V_{\Hcal}$. 

For this purpose we introduce an algebra of meromorphic $\G$-automorphic 
forms: for any integer $d$, we denote by $A_{\Hcal,d}^\G$ the space of
$\G$-invariant  holomorphic functions homogenenous of degree $-d$ 
on $V_{\Hcal}$ that are meromorphic on $V_+$ and have a pole of order at most 
$d$ along the hyperplane sections $H_+$, $H\in\Hcal$. It is clear that 
$A_{\Hcal,\pt}^\G$ is an algebra of holomorphic functions on $\G\bs V_{\Hcal}$
which contains $A_\pt^\G$ as a subalgebra. In particular, this algebra
separates the points of $\G\bs V_{\Hcal}$. We now quote from \cite{looij1}:

\begin{theorem}\label{thm:arrc}
The graded $\CC$-algebra $A_{\Hcal,\pt}^\G$ has finitely many 
generators of positive degree and the resulting complex-analytic injections
\[
\G\bs V_{\Hcal}\to \spec (A_{\Hcal,\pt}^\G),\quad \G\bs \PP (V_{\Hcal})\to \proj (A_{\Hcal,\pt}^\G)
\]
are open embeddings onto Zariski open-dense subsets so that
$\G\bs V_{\Hcal}$ resp.\ $\G\bs \PP (V_{\Hcal})$ acquires the structure of a
quasi-affine resp.\ quasi-projective complex variety. Moreover,
we have natural topological extensions $V_{\Hcal}\subset V_{\Hcal}^*$ and 
$\PP (V_{\Hcal})\subset\PP (V_{\Hcal}^*)$ as $\G$-spaces (described below)
such that these embeddings extend to homeomorphisms 
\[
\G\bs V_{\Hcal}^*\cong \spec (A_{\Hcal,\pt}^\G),
\quad \G\bs \PP( V_{\Hcal}^*)\cong \proj (A_{\Hcal,\pt}^\G)
\]
when the targets are endowed with their Hausdorff topology.
These topological extensions come as $\G$-equivariant stratified spaces which 
result in partitions  of  $\spec (A_{\Hcal,\pt}^\G)$ and  
$\proj (A_{\Hcal,\pt}^\G)$ into subvarieties.
\end{theorem}

\subsection*{Small modification of the Baily-Borel  extension}
The $\G$-extensions in the above theorem share a number of properties with the
Baily-Borel and the toric extensions: the boundary material that we add
comes partitioned into subvarieties (`strata'), where each stratum is given 
as a topological quotient of the extended space. Before we describe the 
extension in question, we first discuss a small modification of the Baily-Borel  extension.

Given $I\in\Ical$, we denote by $\Hcal_{I}$  the collection
of $H\in\Hcal$ containing the isotropic line $I$ and by $I^\Hcal$ the intersection
of $I^\perp$ with all the members of $\Hcal_{I}$. So $I\subset I^\Hcal \subset I^\perp$.
Our rationality assumptions imply that $\G_I$ has finitely many orbits in $\Hcal_I$.
The small modication of the Baily-Borel extension 
amounts to replacing in its definition $I^\perp$ by $I^\Hcal$:
\[
\PP(V_+^{\Hcal}):=\PP(V_+)\sqcup \coprod_{I\in\Ical} \PP(V^+/I^\Hcal)
\]
This extension  comes with a  natural $\G$-invariant topology 
which makes the evident map $\PP(V_+^{\Hcal})\to \PP(V_+^*)$
continuous and its $\G$-orbit space comes with the structure of 
a normal complex-analytic space so that the induced map 
$\G\bs \PP(V_+^{\Hcal})\to \G\bs \PP(V_+^*)$
is an analytic morphism. That morphism is in fact a modification:
it is proper and an isomorphism over the open-dense stratum. The lemma below describes the exceptional fibers and also helps us to understand the behavior of the collection $\Hcal_I$ near $I$.

\begin{lemma}\label{lemma:isotropic} 
The collection $\Hcal_{I}$ induces on $I^\perp$ a finite arrangement (which we denote $\Hcal_I|I^\perp$). 
Any proper intersection $L\not=I^\perp$ of members of $\Hcal_I|I^\perp$ is 
also an intersection of members of $\Hcal_I$ only; moreover, there exists 
an intersection $\tilde L$ of members of $\Hcal_I$ which meets $V_+$ and 
is such that $L=I^\perp\cap \tilde L$.

The fiber of $\G\bs V_+^{\Hcal}\to \G\bs V_+^*$ over a point of the $\CC^\times$-orbit
defined by $I\in\Ical$ is the quotient of an affine space over 
$I^\perp/I^\Hcal$ by a crystallographic group.
\end{lemma}
\begin{proof}
We observe $V_+/I_\Hcal=V/I^\Hcal -I^\perp/I^\Hcal$
is the complement of a linear hyperplane and so its projectivization 
$\PP (V_+/I^\Hcal)=\PP(V/I^\Hcal)-\PP (I^\perp/I^\Hcal)$ is in a natural manner an affine 
space. Let us denote this affine space simply by $A$. 
The group $\G_I$ acts on $A$ through a complex crystallographic group (and so the orbit space of that action is a compact complex-analytic variety (in fact, a quotient of an abelian variety by a  finite group). The exceptional fibers over
the image of $V_+/I^\perp$ are all isomorphic to $\G_I\bs A$ and so the last assertion
of the lemma follows.

Every $H\in \Hcal_I$ determines a hyperplane $A_H$ of $A$ and this defines 
a bijection between $\Hcal_I$ and a collection of affine 
hyperplanes of $A$, a collection we shall denote by $\Hcal |A$. By 
construction, the common intersection of the collection $\Hcal |A$ is empty. 
Notice that for $H_1, H_2\in \Hcal_I$, we have $H_1\cap I^\perp=H_2\cap I^\perp$ 
if and only if $A_{H_1}$ and $A_{H_2}$ are parallel. 

The group $\G_I$  has finitely many orbits in $\Hcal_I$ and hence also in 
$\Hcal_I|I^\perp$ and $\Hcal_I|A$. The assertions of the lemma now follow easily.
The finiteness of $\Hcal_I|I^\perp$ is a consequence of the fact that $\G_I$ acts in 
$I^\perp/I$ through a finite group. This is equivalent to $\Hcal |A$ decomposing into 
finitely many equivalence classes for the relation of parallelism.
If $L\not=I^\perp$ is as in the lemma, then the  collection $\Hcal_L$ of 
$H\in\Hcal$ containing $L$  corresponds to a nonempty finite union of
such equivalence classes; a minimal nonempty intersection of these is an affine
subspace $B$ of $A$ which has $L$ as translation space; the common
intersection $\tilde L$ of the corresponding subset of $\Hcal_L$ will
have the property that $\tilde L$ meets $V_+$ and $\tilde L\cap I^\perp=L$. Finally,
if $B'$ is another such  minimal nonempty intersection distinct from $B$ so that
$B\cap B'=\emptyset$, then the collection of $H\in\Hcal_L$ with $A_H\supset B$ or
$A_H\supset B'$ has $L$ as its common intersection. 
\end{proof}

\begin{remark}\label{remark:cartier}
The modification  $\G\bs \PP(V_+^\Hcal)\to  \G\bs \PP(V^*_+)$ has a simple algebro-geometric meaning: for every $H\in\Hcal$, the evident map $\G_H\bs\PP( H_+)\to \G\bs  \PP(V_+)$ is a finite morphism which extends to the Baily-Borel compactifications:
$\G_H\bs  \PP(H^*_+)\to \G\bs  \PP(V^*_+)$. The latter is also finite with image a hypersurface (in fact, it is a normalization of this image), but that hypersurface
need not support a $\QQ$-Cartier divisor. According to Lemma 5.2 of  \cite{looij1},  the morphism $\G_H\bs  \PP(H^*_+)\to \G\bs  \PP(V^*_+)$ lifts naturally to a morphism
$\G_H\bs  \PP(H^*_+)\to \G\bs  \PP(V^\Hcal_+)$ whose image \emph{does} support a $\QQ$-Cartier divisor. The modification $\G\bs \PP(V_+^\Hcal)\to \G\bs \PP(V^*_+)$  is the smallest one for which the images of all the  $\G_H\bs  \PP(H_+)$ extend to hypersurfaces which are $\QQ$-Cartier  divisors.
\end{remark}

\subsection*{Baily-Borel extension of an arithmetic arrangement complement}
 Let $\Lcal_+(\Hcal)$ denote the collection of all the linear subspaces  
that arise as an intersection of members of $\Hcal$ and meet $V_+$ 
(this includes $V$ as the empty intersection) and denote by 
$\Lcal_+(\Hcal,\Ical)$ the set 
of pairs $(L,I)\in\Lcal_+(\Hcal)\times\Ical$ with 
$L\supset I$. Then the extensions appearing in Theorem \ref{thm:arrc} are
\begin{align*}
V_{\Hcal}^* &=V_{\Hcal}\sqcup 
\left(\coprod_{(L,I)\in\Lcal_+(\Hcal,\Ical)} V_\Hcal/(L\cap I^\perp)\right)
\sqcup \coprod_{L\in\Lcal_+(\Hcal)} (V_\Hcal/L),\\
\PP (V_{\Hcal}^*) &=\PP(V_{\Hcal}) \sqcup 
\left(\coprod_{(L,I)\in\Lcal_+(\Hcal,\Ical)} \PP(V_\Hcal/(L\cap I^\perp))\right)
\sqcup\coprod_{V\not=L\in\Lcal_+(\Hcal)}
\PP ( V_\Hcal/L), 
\end{align*}
with a topology enjoying similar properties as in the Baily-Borel case. Notice,
incidentally, that we recover the latter if $\Hcal$ is empty. 

\subsection*{Structure of the strata}
A valuable piece of information contained in Theorem \ref{thm:arrc} is the
stratified structure it exhibits in the boundary $\proj (A^\G_{\Hcal,\pt})-\G\bs\PP (V_+)$. Any stratum is clearly of the form  $\G_L\bs \PP(V_{\Hcal}/L)$ with $L\in \Lcal_+(\Hcal)$
or  $\G_{L,I}\bs \PP(V_{\Hcal}/L\cap I^\perp)$ with  $(L,I)\in \Lcal_+(\Hcal,\Ical)$.

Let us see what we get in the first case.
Denote by $\Hcal_L$ the collection of members of $\Hcal$ which contain $L$ and 
by $(V/L)_{\Hcal}$ the complement of the union of the members of $\Hcal_L$ in $V/L$.
Then $V/L\cong L^\perp$ is negative definite,  $\Hcal_L$  is finite and $V_{\Hcal}/L=(V/L)_{\Hcal}$. So $\PP(V_{\Hcal}/L)=\PP(V/L)_{\Hcal}$ is a projective arrangement complement. The stabilizer $\G_L$ acts on $L^\perp$ through a finite group and 
hence the stratum $\G_L\bs \PP(V_{\Hcal}/L)$ of $\proj (A^\G_{\Hcal,\pt})$ is well-understood. In particular, its codimension is the dimension of $L$.

In the second case, we only note that the maximal stratum associated to $I$ is the $\G_I$-orbit space of the affine arrangement complement $A_\Hcal:= A-\cup_ {H\in \Hcal_I} A_H$, where $A$ is the affine space and $A_H$ the affine hyperplane that we encountered in the proof of Lemma \ref{lemma:isotropic}. (So this is an open subset
of a finite quotient of an abelian variety). This discussion has an interesting corollary:

\begin{corollary}\label{cor:merom}
Suppose that $\Hcal\not=\emptyset$. Then the codimension of the boundary 
that $\G\bs V_\Hcal$ has in $\proj (A^\G_{\Hcal,\pt})$ is the minimal dimension of a
member of $\Lcal_+(\Hcal)\cup \{ I^\Hcal \}_{I\in\Ical}$ .
In particular, if every one-dimensional intersection of members of $\Hcal$
is  negative definite, then  the boundary of $\G\bs\PP (V_+)$ in 
$\proj (A^\G_{\Hcal,\pt})$ is of codimension $>1$ and the meromorphicity requirement in the definition of $A_{\Hcal,\pt}^\G$ is superfluous: $A_{\Hcal,\pt}^\G$ is simply the algebra of $\G$-invariant holomorphic functions on $V_{\Hcal}$.
\end{corollary}
\begin{proof}
The first statement follows immediately from Theorem 
\ref{thm:arrc} and Lemma \ref{lemma:isotropic}. 

Under the assumption that every one-dimensional intersection of members of $\Hcal$
is  negative definite, the boundary of our compactification is of codimension
$\ge 2$. A meromorphic $\G$-invariant holomorphic functions on $V_{\Hcal}$ which 
is homogeneous of degree $-d$ defines a holomorphic function on 
$\spec (A^\G_{\Hcal,\pt})-\G\bs V_+$. Since $\spec (A^\G_{\Hcal,\pt})$ is normal,
such a function extends to all of $\spec (A^\G_{\Hcal,\pt})$. The homegeneity
assumption implies that this extension is regular, i.e., lies in $A^\G_{\Hcal,d}$.
\end{proof}

\subsection*{Comparison of two compactifications}
There is in general no
natural map $\PP (V_{\Hcal}^*)\to \PP (V^*_+)$ and hence no natural 
morphism $\proj (A^\G_{\Hcal,\pt})\to \proj (A^\G_{\pt})$. As both 
$\proj (A^\G_{\Hcal,\pt})$ and $\proj (A^\G_{\pt})$ are projective compactifications
of the same variety, we can only say that we have a birational map between them.
In \cite{looij1} we explicitly described the graph of the birational map
$\proj (A^\G_{\Hcal,\pt})\dasharrow \proj (A^\G_{\pt})$ and for
its topological  $\G$-equivariant counterpart 
$\PP (V_{\Hcal}^*)\dasharrow  \PP (V^*_+)$. 
We shall not recall this, but it is worth mentioning  an interesting  special case, which is relevant for the example that we will discuss below.

\begin{proposition}\label{prop:contraction}
If  $\Hcal$ has the property that any two distinct members of $\Hcal$ do no intersect each other in $V_+$, then
\begin{enumerate}
\item[(i)] we have a natural map 
$\pi: \PP(V_+^\Hcal)\to \PP(V_\Hcal^*)$ which is continuous $\G$-equivariant and such that the resulting  map  $\G\bs \PP(V_+^\Hcal)\to \G\bs \PP(V_\Hcal^*)$ is a morphism, 
\item[(ii)] the dimension of the boundaries 
$\G\bs \PP(V_+^\Hcal)-\G\bs \PP(V_+)$ 
and $\G\bs \PP(V_\Hcal^*)-\G\bs \PP(V_\Hcal)$  are at most one,
\item[(iii)] the images of the natural morphisms
$\G_H\bs \PP(H_+^*)\to\G\bs\PP(V_+^\Hcal)$ are disjoint if we let $H$ run over a system of representatives  of the $\G$-orbits  in $\Hcal$, and $\pi$ contracts each of these
images.
\end{enumerate}
In particular,  $\G\bs\PP(V_+^\Hcal)$ appears as the normalization of the graph of 
the map of orbit spaces
$\G\bs \PP(V^*_+)\dasharrow \G\bs \PP(V^*_\Hcal)$.
\end{proposition} 

This is a special case of Theorem 5.7 of \cite{looij1}. 
We confine ourselves here to describing the map $\pi $. 
Under the hypothesis of the proposition, the boundary strata of $\PP(V^*_\Hcal)$ consist of  $\PP(V_\Hcal /H)$ (a singletons), $H\in\Hcal$
and $\PP(V_\Hcal  /I^\Ical)$ (a singleton in case $\Hcal_I$ is empty and
an affine line minus a discrete set otherwise), $I\in\Ical$. Now $\pi$ sends $\PP(H_+)$ to
the singleton $\PP(V_\Hcal /H)$; if $I\in\Ical$, then it is identity on
$\PP(V_\Hcal/I^\Hcal)$, and if $H\in\Hcal_I$, then the image of 
$\PP (H_+)$ in $\PP (V_+/I^\Hcal)$ (a singleton)  maps to the singleton $\PP(V_\Hcal /H)$.

\section{The moduli space  plane quartics}
We illustrate the preceding (and especially Proposition \ref{prop:contraction}) with the case of quartic plane curves.

\subsection*{Geometric invariant theory of quartics}
Fix a complex vector space $W$ of dimension three. In what follows
a central role is played by the space of homogeneous quartic forms on $W$, $\sym^4(W^*)$, and so we prefer a briefer name: we abbreviate that space simply by $\Qcal$. The fundamental theorem of geometric invariant theory 
says that each fiber of the natural map  $\Qcal\to \spec \CC[\Qcal]^{\SL (W)}$ contains a unique closed $\SL (W)$-orbit which lies in the closure of all
other orbits in that fiber. In other words, $\spec \CC[\Qcal]^{\SL (W)}$ 
is the orbit space in the separated category, $\SL (W)\bs\bs \Qcal$. The preimage of $0$ is called the \emph{unstable locus} and has the origin of $\Qcal$ as its 
unique closed orbit. The complement of the unstable locus is by definition
the \emph{semistable locus} $\Qcal^\ss$. The 
\emph{stable locus} $\Qcal^\s\subset \Qcal^\ss$ is the set of $F\in \Qcal$ 
for which the orbit map $\SL (W)\to\Qcal$, $g\mapsto gF$, is proper (so that the 
orbit of $F$ is closed in $\Qcal$). Following Mumford these are precisely the 
$F$ which define a reduced quartic curve with only ordinary double
points or (ordinary) cusps.
The closed orbits in the \emph{strictly semistable locus} $\Qcal^\ss-\Qcal^\s$ are the $F\in\Qcal$  of the form 
$(Z_1Z_2-Z_0^2)(sZ_1Z_2-tZ_0^2)$ with $s\not=0$, where $Z_0,Z_1,Z_2$ is a 
coordinate system  for $W$. So
these define quartics that are the union $C'\cup C''$ of two conics with $C'$ smooth,
$C''$ not a double line and for which either $C'$ and $C''$ meet in two points of multiplicity two, or  $C'=C''$. 

The degrees appearing in $\spec \CC[\Qcal]^{\SL (W)}$ are divisible by $3$ because the center $\mu_3$ of $\SL (W)$ acts effectively on $\Qcal$.

We denote by $\Qcal^\circ\subset\Qcal^\s$ the set of $F\in\Qcal$ for which
$C(F)$ is smooth.

\subsection*{Associated K3 surface}
The following construction follows  S.~Kond\=o. 
Fix  $F\in\Qcal^\s$ and denote by $C(F)\subset\PP(W)$ the curve defined by $F$.
The equation $T^4=F$ defines an affine cone $Z(F)$ in $W\oplus\CC$,
whose base at infinity is a quartic surface $S(F)\subset \PP(W\oplus\CC)$.
The projection $S(F)\to \PP(W)$ from $[0:1]\in \PP(W\oplus\CC)$  has as its 
fibers the $\mu_4$-orbits and  $S(F)$ as discriminant. In this way $S(F)$ is a 
$\mu_4$-covering of $\PP(W)$ with total ramification along $C(F)$. Since the 
singularities of $C(F)$ are ordinary double points or
ordinary cusps, those of $S(F)$ are DuVal singularities of type $A_3$ 
(local-analytic equation $w^4=x^2+y^2$) or $E_6$ (local-analytic equation 
$w^4=x^3+y^2$) and hence, if $\tilde S(F)\to S(F)$ resolves these in 
the standard (minimal) manner, then $\tilde S(F)$ is a nonsingular K3 surface. 
This resolution is unique (up to unique isomorphism, to be precise), and so the $\mu_4$-action 
on $S(F)$ lifts to $\tilde S$. The hyperplane class 
$\eta(F)\in H^2(\tilde S(F),\ZZ)$ is a
semipolarization of $\tilde S(F)$ invariant under $\mu_4$ and of degree $4$: 
we have $\eta (F)\cdot \eta(F)=4$ (the dot denotes the intersection form on 
the cohomology). 

We now fix a generator $\mu\in \wedge^3W^*$ and view this 
generator as a translation invariant $3$-form on $W$. Then
the surface $\tilde S(F)$ comes with a holomorphic differential 
defined by
\[
\omega(F):=\res_{S(F)}\left(\frac{\mu}{T^3}\Big|_{Z(F)} \right).
\]
Since $\mu_4$ acts on the last factor $\CC$ by the tautological character,
it acts on the $T$-coordinate by $\bar\chi$ and hence  on $\omega(F)$
by $\bar\chi^{-3}=\bar\chi$.
We often regard $\omega (F)$ as a cohomology class in 
$H^2(\tilde S(F),\CC)_{\bar\chi}$.  Since 
$\omega(F)\wedge\overline{\omega}(F)$ is 
everywhere positive (for the complex orientation), we have
$\omega(F)\cdot \overline{\omega}(F)>0$. 
It is clear that $\tilde S(F)$ is a $K3$-surface. The pull-back 
of the hyperplane class defines a semipolarization $\eta (F)$ of degree four
for $\tilde S(F)$ relative to which the surface is \emph{nonhyperelliptic}: 
$\tilde S(F)$ has no elliptic fibration such that the semipolarization is of 
degree $1$ or $2$ on the fibers. (Otherwise these fibers would map with 
degree $2$ to their image or get contracted.)

If we carry out this construction universally, then we find 
a hypersurface $\Zcal\subset \Qcal^\s\times (W\oplus\CC)$ (defined by the equation 
$F=T^4$), the projection $\Zcal\to\Qcal^\s\times W$ is a $\mu_4$-cover 
which ramifies (totally) over the zero set $\Ccal\subset \Qcal^\s\times W$ of 
$F$ and at infinity we get $\Scal\subset\Qcal^\s\times \PP(W\oplus\CC)$ as a 
$\mu_4$-cover of $\Qcal^\s\times \PP(W)$.
If we let $\GL (W)$ act as the identity on the last factor $\CC$, then
this action preserves $\Zcal$ and its defining equation. We also get a section 
$\omega$ of the relative dualizing sheaf of $\Scal$. This section 
transforms under $\mu_4$ according to the character 
$\det:\GL (W)\to\CC^\times$.

\begin{proposition}\label{prop:assk3}
Let $F\in\Qcal^\s$, denote by $a_1$ resp.\ $a_2$ 
the number of ordinary double points resp.\ ordinary 
cusps of $C$ and put $d:=a_1+2a_2$. Then
 $\mu_4$ acts on $H^2(S(F),\CC)$ with character
$1+(7-d)(\chi +\chi^2+\chi^3)$ and on 
$H^2(\tilde S(F),\CC)$ with character 
\[
1+7(\chi +\chi^2+\chi^3)+ a_1(3-\chi -\chi^2-\chi^3)
+a_2(4-2\chi -2\chi^3).
\]
\end{proposition}
\begin{proof}
The fixed point set of $\mu_4$ in $S=S(F)$ is $C=C(F)$, whereas  the action
of $\mu_4$ on $S-C$ is free. We can now invoke the Lefschetz theorem
which says that the (alternating) character of $\mu_4$ on the total cohomology 
of $S$ takes on $\xi\in \mu_4$ the value $e(S^\xi)$ (the Euler characteristic
of the fixed point set of $\xi$). The latter is $e(S)$ for $\xi =1$ and
$e(C)$ otherwise. In the presence  of $a_1$ ordinary double points and $a_2$ cusps, the Euler characteristic of $C$ resp.\ $S$  is $-4+d$ resp.\ $24-3d$. 
It then easily follows  that the character on the total cohomology is 
$3+(7-d)(\chi +\chi^2+\chi^3)$. So the character on $H^2(S,\CC)$ is as asserted.

We use this to compute the character of $\mu_4$ on $H^2(\tilde S,\CC)$: The difference
between the two is accounted for by  $H^2$ of the exceptional set. 
A double point resp.\ a cusp yields a DuVal curve $D$ of type $A_3$ resp.\ $E_6$.
Its reduced cohomology only lives in dimension $2$ and has as  basis the fundamental classes of the irreducible components of $D$.
A generator of $\mu_4$ leaves in the first case each irreducible component 
invariant (so we get character $3$) and induces in the second case the only non trivial
symmetry of an $E_6$-diagram (so we get character  $4+2\chi^2$). Hence
the character on $H^2(\tilde S,\CC)$ is 
$1+(7-d)(\chi +\chi^2+\chi^3)+a_13+a_2(4+2\chi^2)=
1+7(\chi +\chi^2+\chi^3)+a_1(3-\chi -\chi^2-\chi^3)+ a_2(4-2\chi -2\chi^3)$.
\end{proof}

It is clear that the orbit space $S'(F):=\mu_2\bs S(F)$ is the degree two 
cover of $\PP^2$ which ramifies along $C(F)$.
So $S'(F)$ is  a Del Pezzo surface of degree $2$ which is allowed to have 
DuVal singularities of type $A_1$ and $A_2$.

In Proposition \ref{prop:assk3}, $H^2(S(F),\CC)^{\mu_4}$ is spanned by the 
hyperplane class and $\sqrt{-1}\in\mu_4$ acts on $S(F)$ with Lefschetz number 
$-4+d$. So the following proposition provides a converse in case the latter 
is $\le 0$.

\begin{proposition}\label{prop:converse}
Let $(S,\eta)$ be a polarized K3 
surface  (DuVal singularities allowed) of degree 4 with $\mu_4$-action 
$\rho :\mu_4\to\aut (S,\eta)$ and let $\omega$ be a generator
$\omega$ of the dualizing sheaf of $S$ such that
\begin{enumerate}
\item[(i)] $H^2(S,\CC)^{\mu_4}$ is spanned by $\eta$,
\item[(ii)] the Lefschetz number of $\rho (\sqrt{-1})$ is $\le 0$,
\item[(iii)] $\omega\in H^2(S,\CC)_{\bar\chi}$ and
\item[(iv)] if $\tilde S\to S$ is the minimal resolution, then for every $\eps\in \pic (\tilde S)$
with $\eps\cdot\eps =0$ we have $|\eps\cdot\eta|>2$. 
\end{enumerate}
Then there is an embedding of $S$ in $\PP^3$ such that
$\eta$ is the hyperplane class, the image has an equation
of the form $F(Z_0,Z_1,Z_2)=Z_0^4$, with $F$ defining a stable quartic plane
curve (i.e., $F\in\Qcal^\s$), the $\mu_4$-action on $S$ is the 
restriction of its diagonal action on $\PP^3$ with characters 
$(1,1,1,\chi)$ and $\omega$ is the residue
of $F^{-1}dZ_0\wedge dZ_1\wedge dZ_2$. Moreover the $\SL (3,\CC)$-orbit of $F$ is
unique.
\end{proposition}
\begin{proof}
The last assumption says that $(S,\eta)$ is nonhyperelliptic, i.e., that the degree of every elliptic curve on a minimal resolution $\tilde S$ of $S$ is $\ge 3$. 
According to  Mayer \cite{mayer} and Saint-Donat \cite{stdonat}, the nonhyperellipticity of 
$\tilde S$ implies that the polarization defines up to a projective transformation an embedding $S\subset \PP^3$.
  
Let $S^{\mu_4}\subset S$ denote the fixed point set of the $\mu_4$-action. The
Lefschetz formula says that the Euler characteristic of $S^{\mu_4}$ is the Lefschetz 
number of $\sqrt{-1}\in\mu_4$, which is $\le 0$. Since  $S^{\mu_4}$ is nonempty, it
follows that $S^{\mu_4}$ 
contains an irreducible curve $C$ of positive genus. That curve 
spans a subspace of $\PP^3$ on which $\mu_4$ acts trivially and that subspace is not  
a line. So $\mu_4$ acts trivially on a plane in $\PP^3$. It follows that we can represent 
the $\mu_4$-action on $\CC^4$ by a diagonal action of type 
$(1,1,1,\chi^{\pm 1})$. Then $S$ will have an equation of the form 
$F(Z_0,Z_1,Z_2)=cT^4$ for some $c\in\CC$. Since $S$ has only DuVal 
singularities, we cannot have $c=0$. We rescale the coordinates in such a 
manner that $\omega$ is the residue as above. It is now also clear
that the $\mu_4$-action is of type $(1,1,1,\chi)$.
At the same time we see  that $F$ defines a stable quartic curve.
The uniqueness of the $\SL (3,\CC)$-orbit of $F$ is left to the reader.
\end{proof}

\subsection*{The period map}
We briefly review the period map, again essentially following Kond\=o \cite{kondo}.
Fix an even unimodular lattice $\Lambda$ of signature $(3,19)$ and a
$\eta\in\Lambda$ with $\eta\cdot\eta=4$. As is well-known, such 
a pair $(\Lambda,\eta)$  is unique up to isometry. Consider the collection
of $\mu_4$-actions  $\rho: \mu_4\to \Or (\Lambda)_\eta$ on $\Lambda$ fixing $\eta$
for which
\begin{enumerate}
\item[(i)] $\eta$ spans $\Lambda^{\mu_4}$,
\item[(ii)] the sublattice $\Lambda^{\mu_2}$ is nondegenerate of signature $(1,7)$, and
\item[(iii)] if $\eps\in\Lambda^{\mu_2}$ is such that $\eps\cdot\eps=0$, then 
$\eps\cdot\eta\not= 2$.
\end{enumerate}
Notice that (i) and (ii) imply that $\mu_4$ has character $1+7(\chi+\chi^2+\chi^3)$ on
$\Lambda$. The group $\Or (\Lambda)_\eta$ acts on this set (by composition). 
One can either invoke the surjectivity of the period map and the above discussion 
or use more intrinsically the theory of lattices to see that this action is  transitive. 

We therefore fix one such $\rho$ and we write simply $a'$ for $\rho (\sqrt{-1})a$. 
We let $\Lambda_{\pm}$ denote the set of $a\in\Lambda$ with $a''=\pm a$
(so $\Lambda_+=\Lambda^{\mu_2}$).
According to Kond\=o, $(\Lambda_+,\eta)$ is naturally isomorphic to the Picard
lattice of a Del Pezzo surface of degree two with its anticanonial class, scaled by a factor two: $\Lambda_+$ admits an orthogonal  basis $e_0,\dots, e_7$
with $e_0\cdot e_0=2$, $e_i\cdot e_i=-2$ for $i=1,\dots ,7$ and such that $\eta=3e_0-(e_1+\cdots +e_7)$. 

Notice that $V:=(\Lambda\otimes\CC)_{\bar\chi}$ is a $7$-dimensional
complex vector space. If $u,v\in V$, then $u\cdot v=u'\cdot v'=-u\cdot v$ and so $V$ is totally isotropic  relative to the  $\CC$-bilinear extension of the  form on $\Lambda$.
We have $\bar V= (\Lambda\otimes\CC)_{\chi}$ and $V\oplus \bar V$ is real of signature $(2,12)$. So the  Hermitian extension of the bilinear form
to  $\Lambda\otimes\CC$ has signature $(1,6)$ on $V$. Denote half that 
Hermitian form by $h$ and
let $V_\Ocal$ be the set of $v\in V$ with $\frac{1}{2}(v+\bar v)\in \Lambda_-$. 
So if for instance $\alpha\in \Lambda_-$ is such that $\alpha\cdot \alpha=-2$, then $\alpha\cdot\alpha'=0$, $v:=\alpha + \sqrt{-1}\alpha'\in V_\Ocal$ and we have $h(v,v)=
\frac{1}{2}(\alpha\cdot\alpha +\alpha' \cdot \alpha')=-2$.

Clearly, $V_\Ocal$ is a Hermitian module over the Gaussian integers $\Ocal:=\ZZ[\sqrt{-1}]$.  According to Heckman $V_\Ocal$ admits a $\Ocal$-basis $\{ v_i\}_{i=1}^7$ which at the same time enumerates the vertices of a $E_7$-graph such that
$h(v_i,v_j)$ equals $-2$ when $i=j$,  $0$ when $v_i$ and $v_j$ are not 
connected and $1+\sign(j-i)\sqrt{-1}$ if $v_i$ and $v_j$ are connected. In particular,
the pair $(V,h)$ is naturally defined over $\QQ
(\sqrt{-1})$. Denote by $\tilde \G\subset \Or (\Lambda)$ the group of $\mu_4$-automorphisms of $(\Lambda,\eta)$ and by $\G$ its image in  
the unitary group of $(\Lambda,\eta)$. Both groups are arithmetic and contain
 $\mu_4$ as their center (which acts on $V$ as  group of scalars with 
 character $\bar\chi$). This implies that $\G$ splits as a direct product 
 $\mu_4\times\G_1$: since $V_\Ocal$ is a free $\Ocal$-module of rank $7$, the kernel $\G_1$ of the determinant homomorphism  $\det_\Ocal: \G \to \Ocal^\times=\mu_4$ supplements $\mu_4\subset\G$.
 
Given $F\in\Qcal^\circ$, then a choice of an equivariant isometry
$H^2(S(F);\ZZ)\cong \Lambda$ which takes $\eta (F)$ to $\eta$
will take $\omega (F)$ to a point of $V_+$. The latter's $\G$-orbit 
does not change if we pick another such  isometry and hence we get a well-defined
element of $\G\bs V_+$. The resulting map 
\[
P: \Qcal^\circ\to \G\bs V_+
\]
is analytic and 
constant on the $\SL (W)$-orbits; in fact, for $g\in\GL (W)$, we have 
$P(gF)=\det(g)P(F)$. 
Since $\lambda\in\CC^\times\subset\GL (W)$ takes $F$ to $\lambda^{-4}F$ and 
$P(F)$ to $\lambda^3P(F)$, it follows that $P$ is homogenenous of degree 
$-3/4$ relative to scalar multiplication in $\Qcal$ and $V_+$. 
(There is no contradiction here: if we descend
the $\CC^\times$-action on $V_+$ given by scalar 
multiplication to $\G\bs V_+$, then it is no longer faithful:
the kernel is $\G\cap\CC^\times=\mu_4$.) 

The map $P$ extends across $\PP(\Qcal^\s)$ (for then $S(F)$ only acquires DuVal singularities) and 
yields an analytic map $\PP(\Qcal^\s)\to \G\bs \PP(V_+)$.
The Torelli theorem for $K3$-surfaces implies that this map  
is an open embedding. But it fails to be surjective: if $H$ is a hyperplane of $V$ of signature $(1,5)$ that is orthogonal to some $\eps\in\Lambda$ with $\eps\cdot\eps =0$ and $\eps\cdot \eta=2$, then $\PP(H\cap V_+)$ parametrizes hyperelliptic $K3$-surfaces
and hence its image in $\G\bs \PP(V_+)$ is disjoint with the image of the 
above period map. The following lemma describes the situation in a more precise manner:

\begin{lemma}\label{lemma:dec}
Let $\eps\in \Lambda$ be such that $\eps\cdot\eps=0$, $\eps\cdot\eta=2$ and the span
of $\eta$ and the $\mu_4$-orbit of $\eps$ is of hyperbolic signature.
Then we have the following eigenspace decomposition in $\Lambda\otimes\CC$:
\[
2\eps =\eta+\alpha+\beta=\eta + \tfrac{1}{2}(\alpha -\sqrt{-1}\alpha')+\tfrac{1}{2}(\alpha +\sqrt{-1}\alpha')+\beta, 
\]
with $\alpha,\beta \in \Lambda$ such that $\alpha''=-\alpha$, $\beta'=-\beta$ and
$\alpha\cdot\alpha =\beta\cdot\beta=-2$. 
In particular, $\eps$ is primitive and
the orthogonal projection of $4\eps$ in $V$, $v:=\alpha + \sqrt{-1}\alpha'$, 
lies in $V_\Ocal$ and  satisfies $h(v,v)=-2$. Moreover, $\Ocal v$ is an orthogonal direct summand of $V_\Ocal$.
\end{lemma}

We first prove:

\begin{lemma}\label{lemma:rigid}
Let $L\subset \Lambda$ be a $\mu_4$-invariant sublattice of hyperbolic signature
containing $\eta$ and two distinct isotropic vectors $\eps_1,\eps_2$ with
$\eps_i\cdot\eta=2$. Then $\eps_1\cdot\eps_2=1$. In particular, each $\eps_i$
is primitive.
\end{lemma}
\begin{proof}
Since $\mu_4$ fixes $\eta$,  it will also preserve each connected component 
of $\{ x\in L\otimes\RR -\{0 \} : x\cdot x\ge 0\}$. This implies that 
$\eps_1\cdot\eps_2>0$.  
Now $a_i:=\eta-2\eps_i$ is perpendicular to $\eta$ and we have 
$a_i\cdot a_i=-4$ and $a_1\cdot a_2=-4+4\eps_1\cdot\eps_2$. Since
$a_1$ and $a_2$ span a negative definite lattice, it follows that 
$\eps_1\cdot\eps_2=1$.
\end{proof}

\begin{proof}[Proof of Lemma \ref{lemma:dec}]
We put $\alpha:=\eps-\eps''$ and $\beta:=\eps+\eps''-\eta$ and  so that 
we have the orthogonal decomposition $2\eps =\alpha+\beta+\eta$ 
with $\alpha''=-\alpha$ and $\beta'=-\beta$. If we apply the preceding 
lemma to the pairs $\eps, \eps'$ and $\eps,\eps''$, we find that 
$\alpha\cdot\alpha=\beta\cdot\beta=-2$. Clearly $\alpha +\sqrt{-1}\alpha'$ 
is an eigenvector of our automorphism with eigenvalue $-\sqrt{-1}$.

To prove the last assertion, let $z\in \Lambda_-$. Then we have 
$\alpha \cdot z=(2\eps-\eta-\beta)\cdot z=2\eps\cdot z\in 2\ZZ$ and likewise
$\alpha' \cdot z\in 2\ZZ$. Hence $z$ can be written as $z_1+\frac{1}{2}(\alpha\cdot z)\alpha +\frac{1}{2}(\alpha'\cdot z)\alpha'$ with $z_1\in\Lambda_-$ perpendicular to $\alpha$ and $\alpha'$.
This proves that $\ZZ\alpha+\ZZ\alpha'$ is an orthogonal direct summand of
$\Lambda_-$. Hence $\Ocal v$ is an orthogonal direct summand of
$V_\Ocal$.
\end{proof}

It is clear that in this situation the orthogonal complement $H$ of $v$ in $V$ is also the 
intersection of $V$ with the orthogonal complement of $\eps$ in $\Lambda_\CC$.
Lemma \ref{lemma:dec} tells us that $V_\Ocal$ is the orthogonal direct sum of $H\cap V_\Ocal$ 
and $\Ocal v$, in particular, $H$ is defined over $\QQ (\sqrt{-1})$. 
Kond\=o shows that the underlying $\ZZ$-lattice, the orthogonal complement of $\ZZ \alpha+\ZZ\alpha'$,
is isometric to $U(2)\perp U(2)\perp D_8(-1)$, where $U$ denotes the hyperbolic plane,
$D_8$ the root lattice of type $D_8$, and the number between parenthesis indicates
the scaling factor of the form.
The signature of $H$ is $(1,5)$ and so $H$ is $\G$-rational. We shall refer to a $\eps$ as in this lemma as a  \emph{hyperelliptic vector} and we call the associated $H\subset V$ a  \emph{hyperelliptic hyperplane}. 
We denote the collection of the latter by $\Hcal_h$. This is clearly 
a $\G$-arrangement. In particular, we have defined $V^*_{\Hcal_h}$. 
We verify that this arrangement satisfies the hypotheses of Proposition
\ref{prop:contraction}.

\begin{lemma}\label{lemma:hstrata}
Two distinct members of $\Hcal_h$ do not
meet inside $V_+$. 
\end{lemma}
\begin{proof}
Suppose we have two distinct  hyperelliptic hyperplanes $H_1,H_2$ which
meet inside $V_+$. This means that the corresponding elliptic vectors
$\eps_1,\eps_2$ will be contained in a $\mu_4$-invariant sublattice
$L\subset\Lambda$ of hyperbolic signature. According to Lemma \ref{lemma:rigid}we 
then have $\eps_1^{(k)}\cdot\eps_2^{(l)}=1$ for all $k,l$. It follows that 
$\alpha_i=\eps_i-\eps''_i$ satisfy $\alpha_1^{(k)}\cdot\alpha_2^{(l)}=0$
so that by  Lemma \ref{lemma:dec} we have an orthogonal 
$\mu_4$-invariant decomposition 
\[
\Lambda_-=(\ZZ\alpha_1\perp\ZZ\alpha'_1)\perp (\ZZ\alpha_2\perp\ZZ\alpha'_2
)\perp K
\]
This implies that $M:=U(2)\perp U(2)\perp D_8(-1)$ has an orthogonal direct
summand (of rank one) spanned by a vector $\alpha$ with $\alpha\cdot\alpha=-2$. 
This, in turn, implies that the discriminant 
quadratic form of $M$, $M^*/M\to \QQ/\ZZ$ (the reduction of  
$x\in M^*\mapsto \frac{1}{2}x\cdot x$), represents $-\frac{1}{4}\in \QQ/\ZZ$ (namely its value on $\frac{1}{2}\alpha$). But a straightforward calculation shows that 
this is not the case. This proves the first assertion. 
\end{proof}

The following statements are proved by Kond\=o or are implicit in his
discussion \cite{kondo} (a more detailed discussion can be found in the 
thesis by M.~Artebani \cite{art}):
\begin{enumerate}
\item[(i)] $\G$ acts transitively on
the collection $\Ical$ of degenerate hyperplanes defined over 
$\QQ (\sqrt{-1})$ and
\item[(ii)] $\tilde\G$ acts transitively on the collection of hyperelliptic vectors
and hence $\G$ acts transitively  on the collection $\Hcal_h$ of 
hyperelliptic hyperplanes.
\end{enumerate}
According to Corollary \ref{cor:merom}, 
the algebra of $\G$-invariant holomorphic functions
on $V_{\Hcal_h}$ is  a $\CC$-graded algebra admitting a finite set 
of homogeneous generators of positive degree.
Since $\mu_4$ acts faithfully as a group of scalars on $V$,
$A_{\Hcal_h ,k}^\G$ is zero when $k$ is not a multiple of $4$.
It follows from the two properties above that the projective compactification 
$\PP(\G\bs V_{\Hcal_h})\subset \PP(\G \bs V_{\Hcal_h}^*)=\proj (A_{{\Hcal_h}}^\G)$ 
adds to $\G\bs \PP(V_{\Hcal_h})$ just two strata: a singleton 
$\{\infty_h\}$ (corresponding to a member of $\Hcal_h$) and an affine 
curve $C_h$ (corresponding to some $I^{\Hcal_h}$). The closure of
$C_h$ is the union of these two: $\overline{C_h}=C_h\cup\{ \infty_h\}$.

\begin{theorem} 
Let $\mu_4$ acts on $\Qcal$ by scalar multiplication. Then the period map defines an isomorphism $(\mu_4\times\SL (W))\bs  \Qcal^\s\to \G\bs V_{\Hcal_h}$ and induces an isomorphism of $\CC$-algebras 
\[
\begin{CD}
A_{{\Hcal_h}}^{\G} @>{\cong}>> \CC[\Qcal ]^{\mu_4\times\SL (W)}
\end{CD}
\]
which multiplies degrees by $3$  and gives rise to an 
isomorphism $(\mu_4\times\SL (W))\bs\bs \Qcal\cong \G\bs V_{\Hcal_h}^*$. The points
$\infty_h$  corresponds to the strictly 
semistable orbit defined by $(Z_1Z_2-Z_0^2)^2$ (double conic)  and the curve $C_h$ corresponds to the curve of strictly semistable orbits defined by
$(Z_1Z_2-Z_0^2)(Z_1Z_2-tZ_0^2)$ with $t\not=1,\infty$.
\end{theorem}  
\begin{proof}
The Torelli theorem for $K3$-surfaces implies that 
$\SL (W)\bs \Qcal^\s\to \G\bs V_{\Hcal_h}$ is an isomorphism after 
projectivization. We noticed that it is homogeneous of degree $-3$
relative to the $\CC^\times$-actions on $W$ and $V_+$.
In $\G\bs V_+$ the $\CC^\times$-action is not effective, but has kernel
$\CC^\times\cap \G=\mu_4$. It follows that this map drops to an isomorphism 
$(\mu_4\times\SL (W))\bs  \Qcal^\s\to \G\bs V_{\Hcal_h}$.

Since $\Qcal-\Qcal^s$ is of codimension $>1$ in  $\Qcal$,
$\CC[\Qcal ]^{\mu_4\times\SL (W)}$ is the algebra of regular functions on 
$(\mu_4\times\SL (W))\bs \Qcal^\s$.
Similarly, since $\G\bs V^*_{\Hcal_h}- \G\bs V_{\Hcal_h}$ is of codimension 
$>1$ in $\G\bs V^*_{\Hcal_h}$, $A_{\Hcal_h}^{\G}$ is the algebra of regular 
functions on $\G\bs V_{\Hcal_h}$. It follows that the isomorphism
$(\mu_4\times\SL (W)\bs \Qcal^\s\cong \G\bs V_{\Hcal_h}$
induces an isomorphism $A_{\Hcal_h}^{\G}\cong\CC[\Qcal ]^{\mu_4\times\SL (W)}$
which multiplies the degrees by $3$.
\end{proof}

\begin{question}
It seems likely that the same statements hold for $\SL(W)$
in relation to $\G_1$, so that for instance 
$A_{\Hcal_h}^{\G_1}$ gets identified with $\CC[\Qcal ]^{\SL (W)}$. 
This raises the following
question: given $F\in \Qcal^\circ$, does the naturally defined $\Ocal$-lattice
inside  $H^2(S(F);\CC)_{\bar\chi}$ (of rank 7) have a canonical generator for
its top exterior power?
\end{question}

\begin{remark}
It is clear that the Baily-Borel compactification  $\G\bs\PP (V_+^*)$ of
 $\G\bs\PP (V_+)$  has a unique cusp and so is a one-point compactification.
According to Proposition \ref{prop:contraction}, the small modification 
$\G\bs \PP(V^{\Hcal_h})$  replaces the cusp by a curve and there is 
a natural contraction $\G\bs \PP(V^{\Hcal_h})\to\G\bs \PP(V^*_{\Hcal_h})$ whose 
exceptional divisor is the image of a natural morphism $\G_H\bs \PP(H^*_+)\to 
\G\bs \PP(V^{\Hcal_h})$ (for any $H\in\Hcal_h$). The intersection of this 
exceptional divisor with  $\G\bs \PP (V_+)$ parametrizes the  hyperelliptic 
curves of genus three. 

The small blow-up has a counterpart in the geometric invariant  theory of quartic curves:  Artebani \cite{art} shows in her thesis 
that it is a GIT quotient of a blow-up of the orbit of the double conic in $\PP(\Qcal^\ss)$ and that we thus obtain a compactification of the moduli space of genus three curves.
\end{remark}

\end{document}